\documentclass[11pt,reqno]{article}
\usepackage{amsmath}
\usepackage{amsfonts,amssymb}
\usepackage{amsthm}
\usepackage{graphics,graphicx}
\usepackage{hyperref}
\numberwithin{equation}{section}

\newcommand{\upchi}{\raise1pt\hbox{$\chi$}}

\newcommand{\Kk}{{\mathcal K}}

\newcommand{\jap}[1]{\langle #1 \rangle} 
\newcommand{\norm}[1]{\left\lVert#1\right\rVert}

\newcommand{\Real}{\mathbb R}
\newcommand{\abs}[1]{\left\vert#1\right\vert}
\newcommand{\set}[1]{\left\{#1\right\}}

\newcommand{\grad}{\nabla}

\newcommand{\K}{\mathcal{K}}

\newtheorem{theorem}{Theorem}
\theoremstyle{definition}
\newtheorem{remark}{Remark}

\theoremstyle{lemma}

\newtheorem{proposition}{Proposition}

\theoremstyle{definition}
\newtheorem{definition}{Definition}
\theoremstyle{lemma}

\begin{document}

\title{Large mass global solutions for a class of $L^1$-critical nonlocal aggregation equations and parabolic-elliptic Patlak-Keller-Segel models} 
\author{Jacob Bedrossian\footnote{\textit{jacob@cims.nyu.edu}, Courant Institute of Mathematical Sciences. Partially supported by NSF Postdoctoral Fellowship in Mathematical Sciences, DMS-1103765}}

\date{\today}
\maketitle 

\section*{Abstract}
We consider a class of $L^1$ critical nonlocal aggregation equations with linear or nonlinear porous media-type diffusion which are characterized by a long-range interaction potential that decays faster than the Newtonian potential at infinity.
The fast decay breaks the $L^1$ scaling symmetry and we prove that `sufficiently spread out' initial data, regardless of the mass, result in global spreading solutions. 
This is in contrast to the classical parabolic-elliptic PKS for which essentially all solutions with more than critical mass are known to blow up in finite time. 
In all cases, the long-time asymptotics are given by the self-similar solution to the linear heat equation or by the Barenblatt solutions of the porous media equation. 
 The results with linear diffusion are proved using properties of the Fokker-Planck semi-group whereas the results with nonlinear diffusion are proved using a more interesting bootstrap argument coupling the entropy-entropy dissipation methods of the porous media equation together with higher $L^p$ estimates similar to those used in small-data and local theory for PKS-type equations.  
    
\bigskip

\section{Introduction}
The focus of this work is to study the following general class of equations in $\Real^d$, $d \geq 2$: 
\begin{equation} \label{def:ADD}
\left\{
\begin{array}{l}
  u_t + \grad \cdot (u \grad c) = \Delta u^m, \;\; m \geq 1,\\
  c = \Kk \ast u, \\ 
  u(0) = u_0 \in L_+^1(\Real^d;(1+\abs{x}^2)dx) \;\; d \geq 2, 
\end{array}
\right.
\end{equation}
where  $L_+^1(\Real^d;\mu) := \set{f \in L^1(\Real^d;\mu): f \geq 0}$. In what follows we will always denote $\norm{u(t)}_1 = M$, which is conserved in time for any reasonable notion of solution. 
The prototype for this set of equations is the classical parabolic-elliptic Patlak-Keller-Segel (PKS), which corresponds to the choices $m = 1$ and $\Kk = \mathcal{N}$, where $\mathcal{N}$ denotes the Newtonian potential: 
\begin{equation} \label{def:PKS}
\left\{
\begin{array}{l}
  u_t + \grad \cdot (u \grad c) = \Delta u, \\ 
  -\Delta c = u. 
\end{array}
\right.
\end{equation}
The PKS model is generally considered to be one of the fundamental models of nonlocal aggregation phenomena, especially aggregation via chemotaxis in certain microorganisms  \cite{Patlak,KS,Hortsmann,HandP}. 
Generalizations with nonlinear diffusion (which models an overcrowding effect) and more general nonlocal interactions such as \eqref{def:ADD} have been proposed as models in a variety of ecological systems \cite{Bio,Topaz,Milewski,GurtinMcCamy77}.
Variants of \eqref{def:ADD} and \eqref{def:PKS} also sometimes appear in physical settings \cite{MasoudShelley14,SireChavanis08}; see also \cite{LiebYau87}. 
The class \eqref{def:ADD} is generally characterized by the competition between the tendency for organisms to diffuse (either under Brownian motion when $m = 1$ or to avoid overcrowding when $m > 1$) and the tendency for organisms to aggregate through nonlocal attraction. 
The models can also be seen as the local repulsive limit of inviscid attractive-repulsive aggregation equations which arise both in biology and material science (see e.g. \cite{Burger07,BertozziLaurentRosado10,balague2013nonlocal,balague2013dimensionality} and the references therein).

The wealth of mathematical work on \eqref{def:PKS} and the variants \eqref{def:ADD} is vast so we will not attempt to make a survey here.   
It is well-known that in $\Real^2$, \eqref{def:PKS} is $L^1$ critical (in the sense that the scaling symmetry of \eqref{def:PKS} preserves the $L^1$ norm) and has a critical mass phenomena (see e.g. \cite{Dolbeault04,BlanchetEJDE06}): if $\norm{u}_1 = M < 8\pi$ then the solution is global and converges to the unique, self-similar spreading solution whereas if $\norm{u}_1 = M > 8\pi$ then the solution blows up in finite time (at least if it has a finite second moment). 
Solutions with exactly critical mass exhibit a variety of possible behaviors including infinite-time aggregation \cite{BlanchetCarrilloMasmoudi08} and convergence to stationary solutions \cite{BlanchetCarlenCarrillo10}. 
In $\Real^3$, \eqref{def:PKS} is $L^1$ (and free energy) supercritical and little beyond small $L^{3/2}$ global existence results (see e.g. \cite{Biler95,Nagai95,Corrias04}) and large $L^{3/2}$ finite time blow up results is known ($L^{3/2}$ is the critical Lebesgue space).
In $\Real^d$ for $d \geq 3$, the choice $\Kk = \mathcal{N}$ and $m = 2-2/d$ is $L^1$ critical, and it was shown in \cite{Blanchet09} that \eqref{def:ADD} with these choices has properties similar to those of \eqref{def:PKS} in $\Real^2$: there exists a critical mass $M_c$ such that if $\norm{u}_1 = M < M_c$ then the solution is global and converges to self-similar spreading solutions whereas if $\norm{u}_1 = M > M_c$, then at least large classes of solutions are known to blow up in finite time (see \cite{Blanchet09,BedrossianKim10}).
Critical mass phenomena also occurs in the more general class \eqref{def:ADD} for suitable choices of $\K$ and $m$ (including also more general filtration equation-type diffusion) \cite{BRB10,BR11,KarchSuzuki10}.

The purpose of this work is to show that for the $L^1$ critical models ($m = 2-2/d$ in $d \geq 2$), if $\K$ decays faster than the Newtonian potential at infinity (in the sense that $\norm{\grad \Kk}_{q} < \infty$ for some $q < \frac{d}{d-1}$),
 then unlike the scale-invariant case, all sufficiently spread out solutions are global and converge to the self-similar spreading solution of the homogeneous diffusion equation $u_t = \Delta u^{2-2/d}$. 
In particular, this covers the well-known case of parabolic-elliptic PKS with lower order degradation term in $\Real^d$, $d \geq 2$ (which is known to have finite time blow-up solutions for all values of $M > M_c$):
\begin{equation} \label{def:PKSalpha}
\left\{
\begin{array}{l}
  u_t + \grad \cdot (u \grad c) = \Delta u^{2-2/d}, \;\;\;\; d \geq 2, \\ 
  -\Delta c + \alpha c = u, \;\;\; \alpha > 0. 
\end{array}
\right.
\end{equation}
The results and proofs are perturbative in nature, treating \eqref{def:ADD} as a perturbation of the diffusion equation in forward self-similar variables. 
Usually in such perturbative settings, the mass (or size of $\K$) is required small, as in for example \cite{CanizoCarrilloSchonbek,SugiyamaDIE07,LuckhausSugiyama07,BedrossianIA10}.
However, here the small parameter that controls the nonlocal aggregation term is basically a measure of the characteristic length-scale of the initial data relative to $\norm{\grad \Kk}_{q}$ for some $q < \frac{d}{d-1}$ (which serves to measure the strength of the attraction on large length-scales) and some appropriate quantification of the size of the initial data. 
We remark that these results are somewhat analogous to behavior observed in the parabolic-parabolic Keller-Segel models \cite{corrias2014existence,BilerCorriasDolbeault11}, where the characteristic time-scale of chemo-attractant diffusion can be used as the small parameter.

That the long-time asymptotics should be governed only by the diffusion equation as $t \rightarrow \infty$ 
has already been observed in, for example, \cite{LuckhausSugiyama06,LuckhausSugiyama07,BedrossianIA10,CanizoCarrilloSchonbek}.
The present work need only concentrate on extending the range of examples where strong decay estimates are known; indeed, for the cases we will study it was shown in \cite{BedrossianIA10} that sufficiently strong decay estimates imply that the solutions converge to the self-similar spreading solution of the diffusion equation.

We will restrict our attention to interaction potentials $\Kk$ that satisfy basic regularity requirements (this definition is originally from \cite{BRB10}). 
Note that while it is not necessary for this work to require $\K$ to be radially non-increasing, which corresponds to $\K$ being purely attractive, the results are mostly interesting when $\K$ is attractive as this is opposing the diffusion.
\begin{definition}[Admissible Kernel] \label{def:adm}
 We say a kernel $\K \in C^3 \setminus \set{0}$ is \emph{admissible} if $\K \in W^{1,1}_{loc}(\Real^d)$ and the following holds:
\begin{itemize}
\item[\textbf{(KN)}] $\K$ is radially symmetric, $\K(x) = k(\abs{x})$ and $k(\abs{x})$ is monotone in a neighborhood of $x= 0$. 
\item[\textbf{(MN)}] $k^{\prime\prime}(r)$ and $k^\prime(r)/r$ are monotone on $r \in (0,\delta)$ for some $\delta > 0$. 
\item[\textbf{(BD)}] $\abs{D^3\K(x)} \lesssim \abs{x}^{-d-1}$. 
\end{itemize}
\end{definition}
The definition ensures that the kernel is radially symmetric, well-behaved at the origin and has second derivatives which define bounded singular integral operators on $L^p$ for $1 < p < \infty$.
It is important to note that all admissible kernels satisfy $\nabla \K \in L^{\frac{d}{d-1},\infty}$, where $L^{p,\infty}$ denotes the weak-$L^p$ space, making the Newtonian potential effectively the most singular of admissible kernels \cite{BRB10}.
Provided $\K$ is admissible, for a given initial condition $u_0(x) \in L_+^1(\Real^d; (1 + \abs{x}^2)dx ) \cap L^\infty(\Real^d)$, 
\eqref{def:ADD} has a unique, local-in-time weak solution which satisfies $u(t) \in C([0,T);L_+^1(\Real^d;(1+\abs{x}^2)dx) )\cap L^\infty((0,T)\times\Real^d)$ for some $T \leq \infty$ (see e.g. \cite{BRB10,BertozziSlepcev10,BlanchetEJDE06,SugiyamaDIE07,BR11}).

In the case of linear diffusion, we will be using strong contractive properties of the Fokker-Planck semi-group which rely on a spectral gap for the associated elliptic problem (see Proposition \ref{prop:Stau}). This generally requires some kind of weighted space; here we define the weighted $L^2$ norm: 
\begin{align*}
\norm{f}_{L^2(\beta)}^2 & = \int \left(1 + \abs{x}^2\right)^{2\beta} \abs{f(x)}^2 dx,  
\end{align*}
with the space $L^2(\beta) = \set{f \in L^1: \norm{f}_{L^2(\beta)} < \infty}$. In what follows denote $\jap{x} = (1+\abs{x}^2)^{1/2}$.  
The statement of Theorem \ref{thm:linear} is given below. 

\begin{theorem}[Linear diffusion] \label{thm:linear}
Let $d = 2$, $m = 1$ and suppose $\Kk$ satisfies Definition \ref{def:adm} and $\norm{\grad \Kk}_q < \infty$ for some $q < 2$.
Then for all $f \in L_+^1 \cap L^2(\beta)$ for some $\beta > 2$, there exists a $\lambda_0 = \lambda_0(\norm{f}_{L^2(\beta)},\norm{f}_1,\beta,\Kk)$ such that if $\lambda > \lambda_0$ and we take the initial data in \eqref{def:ADD} to be
\begin{align} 
u_0(x) = \frac{1}{\lambda^2}f \left(\frac{x}{\lambda}\right), \label{eq:defu0}
\end{align}
then the corresponding solution to \eqref{def:ADD} is global and satisfies the $L^\infty$ decay estimate for $t \geq 1$:
\begin{align} 
\norm{u(t)}_\infty & \lesssim t^{-1}.  \label{ineq:decaylin}
\end{align}
If $\abs{\grad\Kk(x)} \lesssim \abs{x}^{-\gamma}$ for some $\gamma > 1$ then we have
the convergence to self-similarity: for all $\delta > 0$,  
\begin{align}
\norm{u(t) - e^{t\Delta}u_0}_1 & \lesssim_\delta (1+t)^{-\frac{1}{2}\min\left(1,\gamma-1\right) + \delta}. \label{ineq:asymptlin}
\end{align} 
\end{theorem} 

To state our result regarding nonlinear diffusions, recall the self-similar Barenblatt solution of the porous media equation  for $m = 2-2/d$ \cite{VazquezPME}: 
\begin{equation}
\mathcal{U}(t,x;M) = t^{-1}\left( C_1 - \frac{(m-1)}{2md}\left(\frac{\abs{x}}{t^{\frac{1}{d}}}\right)^2\right)_+^{\frac{1}{m-1}}, \label{def:UnonLin}
\end{equation}
where $C_1$ is determined from the conservation of mass. 
Then our result on nonlinear diffusion is stated below. 
The proof is a bootstrap argument that couples a high $L^p$ estimate of the type that arises in the perturbative local or small-data data theory of \eqref{def:ADD} (see e.g. \cite{JagerLuckhaus92,Kowalczyk05,CalvezCarrillo06,Blanchet09,Corrias04,BedrossianIA10}) together with an entropy-entropy dissipation argument based on the inequalities for the porous media equation (see e.g. \cite{CarrilloToscani00,CarrilloEntDiss01}), sometimes considered the nonlinear analogue of a spectral gap. 
That \eqref{ineq:LinftyNL} implies \eqref{ineq:BBlattNL} a posteriori is proved in \cite{BedrossianIA10} using entropy methods (see also \cite{CanizoCarrilloSchonbek}), however, the proof of Theorem \ref{thm:nonlinear} is the only 
example, to the author's knowledge, of a method for PKS-type equations that couples the entropy methods together with perturbative higher $L^p$ estimates to prove a decay estimate of the type \eqref{ineq:LinftyNL}. 

\begin{theorem}[Nonlinear diffusion] \label{thm:nonlinear}
Let $d \geq 3$, $m = 2-2/d$ and suppose $\Kk$ satisfies Definition \ref{def:adm} and $\norm{\grad \Kk}_q < \infty$ for some $q < \frac{d}{d-1}$. 
Then for all $f \in L_+^1(\Real^d; (1+\abs{x}^2)dx) \cap L^\infty$, there exists a $\lambda_0 = \lambda_0(f,\Kk,d)$ such that 
if $\lambda > \lambda_0$ and we take the initial data in \eqref{def:ADD} to be
\begin{align} 
u_0(x) = \frac{1}{\lambda^{d}}f \left(\frac{x}{\lambda}\right), \label{def:u0NL}
\end{align}
then the corresponding solution to \eqref{def:ADD} is global and satisfies the $L^\infty$ decay estimate: 
\begin{align} 
\norm{u(t)}_\infty & \lesssim (1+t)^{-1}.  \label{ineq:LinftyNL}
\end{align}
If $\abs{\grad\Kk(x)} \lesssim \abs{x}^{-\gamma}$ for some $\gamma > d-1$ then we have
the convergence to self-similarity: for all $\delta > 0$,  
\begin{align}
 \norm{u(t) - \mathcal{U}(t,x;M)}_1 & \lesssim_\delta (1+t)^{-\frac{1}{d}\min\left(1,\gamma - d + 1\right) + \delta}. \label{ineq:BBlattNL}
\end{align} 
\end{theorem}

\begin{remark} 
For $L^1$ supercritical cases $1 \leq m < 2-2/d$ (for example the case of parabolic-elliptic PKS in $\Real^3$), 
both Theorems \ref{thm:linear} and \ref{thm:nonlinear} are immediate from small data $L^{\frac{d(2-m)}{2}}$ global existence results even in the case $\Kk = \mathcal{N}$ (see for example \cite{Corrias04,SugiyamaDIE07,SugiyamaDIE06,BedrossianIA10}).  
For more information on supercritical cases, see also \cite{bian2013dynamic} and the references therein. 
\end{remark}

\begin{remark} 
For subcritical problems $m > 2-2/d$ the question of long time behavior has a number of gaps as the aggregation can dominate on large length-scales in these cases. 
To the author's knowledge, no decaying solution for \eqref{def:ADD} with $m > 2-2/d$ has ever been exhibited for an attractive choice of $\Kk$ (e.g. $\grad \Kk \cdot x \leq 0$).  
It is known that in the case $2-2/d < m < 2$, stationary solutions exist for sufficiently large mass for basically all purely attractive choices of $\Kk$ \cite{LionsCC84} (in fact this is true over the entire range $1 < m < 2$ depending on the singularity of the kernel).   
The case $m = 2$ is critical from this perspective \cite{BedrossianGlobalMin10,burger2011stationary}  and in the case $m > 2$ there exists stationary solutions for all values of the mass for basically all radially-symmetric, attractive $\Kk$ \cite{BedrossianGlobalMin10}. 
In some cases, convergence to stationary solutions has been established \cite{KimYao11}.
\end{remark}

\begin{remark} 
If $\gamma \geq d$ then the convergence rates in \eqref{ineq:asymptlin} and \eqref{ineq:BBlattNL} are nearly optimal in the sense that they match the rate of the diffusion equation (up to the $\delta$) \cite{CarrilloToscani00,VazquezPME}. 
In both Theorems \ref{thm:linear} and \ref{thm:nonlinear}, if $\grad \K \in L^1$ we may take $\gamma = d$ in the statement. 
\end{remark}

\begin{remark}
Note that the regularity of $\K$ is essentially irrelevant, it is only the decay at infinity (as long as $\K$ is not more singular than the Newtonian potential). 
For example, both the statements and the proofs of Theorems \ref{thm:linear} or \ref{thm:nonlinear} are the same regardless if we are considering $\K(x) = e^{-\abs{x}^2}$ or $\K$ the fundamental solution of $-\Delta c + \alpha c = 0$ for $\alpha > 0$ and there is no obvious simplification possible in the case of the former. 
\end{remark}

\section{Linear diffusion} \label{sec:lin}
Define the Fokker-Planck operator and linear semi-group
\begin{align*} 
Lf & = \Delta f + \frac{1}{2}\grad \cdot (\xi f) \\ 
S(\tau) & = e^{\tau L}. 
\end{align*} 
We will use some of the following properties of the linear propagator $S(\tau)$ in $L^2(\beta)$, studied in \cite{GallayWayne02}. 
\begin{proposition}[Properties of $S(\tau)$ (see \cite{GallayWayne02})] \label{prop:Stau}
Fix $\beta > 1$. Then, 
\begin{itemize}
\item[(i)] $S(\tau)$ defines a strongly continuous semi-group on $L^2(\beta)$ and for all $w \in L^2(\beta)$, 
\begin{equation}
\norm{S(\tau)w}_{L^2(\beta)} \lesssim \norm{w}_{L^2(\beta)}, \;\;\; \norm{\grad S(\tau) w}_{L^2(\beta)} \lesssim \frac{1}{a(\tau)^{1/2}}\norm{w}_{L^2(
m)},  \label{ineq:SgradHyper}
\end{equation}
for all $\tau > 0$ and where $a(\tau) = 1- e^{-\tau}$. 
\item[(ii)] If $\beta > 2$ and $w \in L_0^2(\beta)$, then
\begin{equation}
\norm{S(\tau)w}_{L^2(\beta)} \lesssim e^{-\tau/2}\norm{w}_{L^2(\beta)}, \;\;\; \forall \tau>0.
\end{equation}
\item[(iii)] If $q \in [1,2]$ then for all $w \in L^q(\beta)$ and $\tau > 0$, 
\begin{align}
\norm{S(\tau)w}_{L^2(\beta)} \lesssim \frac{1}{a(\tau)^{\frac{1}{q} - \frac{1}{2}}}\norm{w}_{L^q(\beta)} \\ 
\norm{\grad S(\tau) w}_{L^2(\beta)} \lesssim \frac{1}{a(\tau)^{\frac{1}{q}}}\norm{w}_{L^q(\beta)}. \label{ineq:gradSDecay}
\end{align}
\end{itemize}
Note that
\begin{equation}
\grad S(\tau) = e^{\tau/2}S(\tau) \grad. \label{eq:SgradCommute} 
\end{equation}
\end{proposition}   

With Proposition \ref{prop:Stau}, we may prove Theorem \ref{thm:linear} with a short perturbation argument.
\begin{proof}[\textbf{(Proof of Theorem \ref{thm:linear})}]  
Denote $u(t,x)$ to be the unique solution to \eqref{def:ADD} with initial data \eqref{eq:defu0}, which is known to exist on some time interval $[0,T_{\max})$ by local well-posedness.  
Define the parameter $T>0$ to be chosen large later: 
\begin{align*} 
T = (\lambda^2-1). 
\end{align*}
Define the self-similar variables $(\tau,\xi)$, 
\begin{align*}
\xi & = ((t+T) + 1)^{-1/2}x \\ 
\tau & = \log \left((t+T) + 1\right),
\end{align*}
together with the rescaled solution 
\begin{equation*}
\theta(\tau,\xi) = ((t+T) + 1) u(t,x),  
\end{equation*}
which is defined on the time interval $[\tau_0, \tau_{\max})$, where 
\begin{align*} 
\tau_0 & = \log(T + 1) \\  
\tau_{\max} & = \log((T_{\max} + T) + 1). 
\end{align*} 
In these variables, \eqref{def:ADD} with initial data \eqref{eq:defu0} becomes the system
\begin{subequations}
  \label{sys:selfsimilarLin2}
  \begin{align}
    \label{eq:selfsimilarLin2}
    & \theta_\tau + \grad \cdot (\theta e^{\tau/2}(\grad \Kk)(e^{\tau/2} \cdot) \ast \theta) = \Delta \theta + \frac{1}{2}\grad \cdot (\xi \theta) \\ 
     &\theta(\tau_0,\xi) = f(\xi). 
    \label{eq:selfsimilarICLin2}
  \end{align}
\end{subequations}
  The idea behind the introduction of $T$ is that if $u_0$ has a characteristic length scale $O(\sqrt{T})$, then $\theta(\tau_0)$ has a characteristic length scale of $O(1)$. The parameter $T$ will eventually be required large to ensure that the initial data lives on a much larger length-scale than the interaction range of the potential. 

Applying Duhamel's formula to \eqref{sys:selfsimilarLin2} gives
\begin{align*}
\theta(\tau) = S(\tau-\tau_0)f - \int_{\tau_0}^\tau S(\tau-s)\left[\grad \cdot (\theta e^{s/2}(\grad \Kk)(e^{s/2} \cdot) \ast \theta(s))\right] ds. 
\end{align*} 
We will be essentially linearizing around the approximate solution $S(\tau-\tau_0)f$.  
Let $[\tau_0,\tau_\star]$ be the largest connected, closed interval such that
\begin{align}  
\norm{\theta(\tau) - S(\tau - \tau_0)f}_{L^2(\beta)} \leq 4, \label{ineq:linBootstrap}
\end{align}
which is well-defined and non-empty by the continuity in time of $\theta(\tau)$ and $S(\tau)$ (Proposition \ref{prop:Stau}). 
Moreover, by standard propagation of regularity, the solution $\theta(\tau)$ is $C^\infty$ for $\tau \in (\tau_0,\tau_\star]$.  
Using the crucial decay estimate \eqref{ineq:gradSDecay}, we deduce 
\begin{align} 
\norm{\theta(\tau) - S(\tau - \tau_0)f}_{L^2(\beta)} & \leq \norm{\int_{\tau_0}^\tau S(\tau-s)\left[\grad \cdot (\theta e^{s/2}(\grad \Kk)(e^{s/2} \cdot) \ast \theta(s)) \right] ds}_{L^2(\beta)} \nonumber \\ 
& \lesssim \int_{\tau_0}^\tau \frac{e^{-\frac{1}{2}(\tau-s)}}{a(\tau-s)^{3/4}}\norm{\theta e^{s/2}(\grad \Kk)(e^{s/2} \cdot) \ast \theta}_{L^{4/3}(\beta)} ds. \label{ineq:Duhamel1}
\end{align} 
By H\"older's inequality: 
\begin{align} 
\norm{\jap{\xi}^{m} \theta e^{s/2}(\grad \Kk)(e^{s/2} \cdot) \ast \theta}_{4/3} \leq \norm{\theta}_{L^{2}(\beta)}\norm{e^{s/2}\grad \Kk(e^{s/2} \cdot) \ast \theta}_{4}. \label{ineq:Holder}
\end{align}
The key here is to use Young's inequality and put $\grad \Kk$ in an $L^z$ space with $z < 2$, breaking the scale invariance that would be present if $\Kk$ were the Newtonian potential (in which case we would only have $\grad \Kk \in L^{2,\infty}$). 
Since $\grad \Kk \in L^{2,\infty}$, by interpolation, $\grad \K$ is in every $L^z$ space with $z \in [q,2)$. 
Therefore, by choosing $q \leq z < 2$ we may ensure by Young's inequality that, 
for some $1 < p < 2$ we have 
\begin{align*} 
\norm{e^{s/2} (\grad \Kk(e^{s/2} \cdot)\ast \theta)}_{4} \lesssim \norm{\theta}_{p} \norm{e^{s/2} \grad \Kk(e^{s/2} \cdot)}_{z} = e^{\frac{s}{2}\left(1 - \frac{2}{z}\right)} \norm{\theta}_{p} \norm{\grad \Kk}_{z}.
\end{align*} 
Since $p < 2$ and $\beta > 2$, by H\"older's inequality we have $\norm{\theta}_p \lesssim_\beta \norm{\theta}_{L^2(\beta)}$, 
so by $\grad \Kk \in L^z$ we have 
\begin{align*} 
\norm{e^{s/2} (\grad \Kk(e^{s/2} \cdot)\ast \theta)}_{4} \lesssim e^{\frac{s}{2}\left(1 - \frac{2}{z}\right)} \norm{\theta}_{L^2(\beta)}. 
\end{align*} 
This exponential decay factor introduces the small parameter we can exploit to close the perturbation argument.  
Using this together with \eqref{ineq:Holder} and \eqref{ineq:Duhamel1} gives us 
\begin{align*} 
\norm{\theta(\tau) - S(\tau - \tau_0)f}_{L^2(\beta)} \lesssim e^{\left(1 - \frac{2}{z}\right)\frac{\tau_0}{2}}\int_{\tau_0}^\tau \frac{e^{-\frac{1}{2}(\tau-s)}}{a(\tau-s)^{3/4}}\norm{\theta(s)}^2_{L^2(\beta)} ds.  
\end{align*} 
Therefore, by the bootstrap hypothesis \eqref{ineq:linBootstrap}, 
\begin{align*} 
\norm{\theta(\tau) - S(\tau - \tau_0)f}_{L^2(\beta)} & \lesssim e^{\left(1 - \frac{2}{z}\right)\frac{\tau_0}{2}} \sup_{s \in (\tau_0,\tau_\star)} \norm{\theta(s)}_{L^2(\beta)}^2 \\
& \lesssim e^{\left(1 - \frac{2}{z}\right)\frac{\tau_0}{2}}\left( 1 + \sup_{s \in (\tau_0,\tau_\star)}\norm{S(\tau - \tau_0)f}_{L^2(\beta)}^2\right). 
\end{align*}
Applying \eqref{ineq:SgradHyper} from Proposition \ref{prop:Stau} implies 
\begin{align*} 
\norm{\theta(\tau) - S(\tau - \tau_0)f}_{L^2(\beta)} & \leq C_1 e^{\left(1 - \frac{2}{z}\right)\frac{\tau_0}{2}} + C_2e^{\left(1 - \frac{2}{z}\right)\frac{\tau_0}{2}} \norm{f}^2_{L^2(\beta)},  
\end{align*} 
where both $C_1$ and $C_2$ are independent of $f$, $\tau_0$ and $\tau_\star$ (they depend only on $\K$, $q$, $\beta$ and the constants coming from Proposition \ref{prop:Stau}).  
By assumption, $\norm{f}_{L^2(\beta)} < \infty$ and hence we may fix $\tau_0$ depending only on the constants $C_i$ and $\norm{f}_{L^2(\beta)}$ such that on $[\tau_0,\tau_\star)$ there holds, 
\begin{align*} 
\norm{\theta(\tau) - S(\tau - \tau_0)f}_{L^2(\beta)} < 2. 
\end{align*}
Therefore, a continuity argument implies that $\tau_\star = \tau_{\max}$ and since $L^2(\beta)$ is a higher $L^p$ space than the critical $L^1$ space, it is standard that the solution is global: $\tau_{\max} = \infty$
and $\norm{\theta(\tau) - S(\tau - \tau_0)f}_{L^2(\beta)} < 2$ for all time. 
The uniform bound in $L^2(\beta)$ on $\theta$ implies the $L^\infty$ decay estimate \eqref{ineq:decaylin} by Theorem 1 (ii) in \cite{BedrossianIA10}, 
and the convergence to self-similarity \eqref{ineq:asymptlin} follows from Theorem 2 or 3 in \cite{BedrossianIA10} (one could alternatively use a second argument via Duhamel's principle as in the methods of \cite{CanizoCarrilloSchonbek}, which might be more natural for linear diffusion).  
\end{proof}

\section{Nonlinear diffusion}   
It is clear that the proof of Theorem \ref{thm:linear} does not apply at all as it depends on the decay estimates of the 
Fokker-Planck semi-group, which are the consequence of an appropriate spectral gap for $L$ in $L^2(\beta)$ (see \cite{GallayWayne02}). 
We instead use the entropy-entropy dissipation inequalities for the porous media equation (see e.g. \cite{CarrilloToscani00,CarrilloEntDiss01}).  
In similarity variables (\cite{VazquezPME,CarrilloToscani00} or \eqref{def:NLsim} below with $T = 0$), the diffusion equation $u_t = \Delta u^{2-2/d}$ is transformed into the nonlinear Fokker-Planck equation: 
\begin{align} 
\theta_\tau = \Delta \theta^{2-2/d} + \grad \cdot (\xi \theta), \label{def:HFkP}
\end{align}
where $\theta(\tau,\xi) = e^{\tau d} u(t,x)$. 
Define the entropy functional 
\begin{equation}
H[\theta] = \frac{1}{m-1}\int \theta^m(\xi) d\xi + \frac{1}{2}\int \abs{\xi}^2 \theta(\xi) d\xi, \label{def:H}
\end{equation}
and the entropy production functional, 
\begin{equation}
I[\theta] = \int \theta\abs{\frac{m}{m-1}\nabla \theta(\xi)^{m-1} + \xi}^2 d\xi. \label{def:I}
\end{equation}
These entropies were originally introduced for studying the porous media equation in \cite{Newman84,Ralston84}.
It is well known that \eqref{def:H} is displacement convex \cite{McCann97} and that \eqref{def:HFkP} is a gradient flow for \eqref{def:H} in the Euclidean Wasserstein distance \cite{Otto01}.
Denote by $\theta_M$ the unique minimizer of the functional \eqref{def:H} with fixed mass $M$ (which is simply the Barenblatt solution \eqref{def:UnonLin} of mass $M$ written in similarity variables) and define the relative entropy
\begin{align*}
H[\theta|\theta_M] = H[\theta] - H[\theta_M] \geq 0. 
\end{align*} 
The functionals are all related by the following: if $\theta(\tau,\xi)$ solves \eqref{def:HFkP}, then 
\begin{equation}
\frac{d}{d\tau}H[\theta(\tau)|\theta_M] = -I[\theta(\tau)]. \label{eq:Hevo}
\end{equation}
Then we have the following, which generalizes the Gross logarithmic Sobolev inequality \cite{Gross75} (see also \cite{DelPinoDolbeault02}). 
\begin{proposition}[Generalized Gross Logarithmic Sobolev Inequality \cite{CarrilloToscani00,CarrilloEntDiss01,DelPinoDolbeault02,Gross75}] \label{thm:rel_entropy}
Let $f\in L_+^1(\Real^d)$ with $\norm{f}_1 = M$. Then, 
\begin{equation}
H[f|\theta_M]  \leq \frac{1}{2} I[f]. \label{ineq:relative_entropy}
\end{equation} 
\end{proposition}
Equations \eqref{eq:Hevo} and \eqref{ineq:relative_entropy}, together with a suitable generalization of the Csiszar-Kullback inequality \cite{Csiszar67,Kullback67,CarrilloToscani00,CarrilloEntDiss01}, provide a sharp quantitative estimate on the rate of convergence of solutions to \eqref{def:HFkP} to $\theta_M$ in $L^1$. Upon transforming back to the original variables, this becomes the convergence to self-similarity for the porous media equation. 

To prove Theorem \ref{thm:nonlinear}, we will begin as in \eqref{eq:Hevo} but will encounter an error term that requires control on a higher $L^p$ norm. 
To control this, we couple the entropy-entropy dissipation argument with the truncated $L^p$ estimate methods which are classical in the study of PKS and its variants. For example, related arguments can be found in \cite{JagerLuckhaus92,Kowalczyk05,CalvezCarrillo06,Blanchet09,SugiyamaDIE07,BedrossianIA10}.
These methods allow to propagate arbitrary $L^p$ estimates provided some uniform equi-integrability is known (see \cite{CalvezCarrillo06}), which here is provided in turn by control on the relative entropy. 
In order to close the bootstrap, the small parameter employed is the length-scale of the initial data.  
  
\begin{proof}[\textbf{(Proof of Theorem \ref{thm:nonlinear})}] 
Denote $u(t,x)$ to be the unique solution to \eqref{def:ADD} with initial data \eqref{def:u0NL}, which is known to exist on some time interval $[0,T_{\max})$ by local well-posedness. 
Define the parameter $T>0$ to be chosen large later: 
\begin{align*} 
T = \frac{1}{d}(\lambda^d-1). 
\end{align*}
As in the beginning of the proof of Theorem \ref{thm:linear}, define the self-similar variables $(\tau,\xi)$ (we remark that the slightly different convention in the definition depending on $d$ holds no real significance): 
\begin{subequations} \label{def:NLsim}
\begin{align}
\xi & = (d(t+T) + 1)^{-1/d}x, \\ 
\tau & = \frac{1}{d}\log\left(d(t+T) + 1\right), \\ 
\theta(\tau,\xi) & = (d(t+T) + 1) u(t,x), 
\end{align}
\end{subequations}
which is defined on the time interval $[\tau_0, \tau_{\max})$, where 
\begin{align*} 
\tau_0 & = \frac{1}{d}\log\left(dT + 1\right), \\
\tau_{\max} & = \frac{1}{d}\log\left(d(T_{\max}+T) + 1\right).
\end{align*} 
Written with \eqref{def:NLsim}, \eqref{def:ADD} with initial data \eqref{def:u0NL} becomes
\begin{subequations}
  \label{sys:selfsimilarNL2}
  \begin{align}
    \label{eq:selfsimilarNL2}
    & \theta_\tau + \grad \cdot (\theta e^{(d-1)\tau}(\grad \Kk)(e^\tau \cdot) \ast \theta) = \Delta \theta^m + \grad \cdot (\xi \theta) \\ 
     &\theta(\tau_0,\xi) = f(\xi). 
    \label{eq:selfsimilarICNL2}
  \end{align}
\end{subequations}
By the regularity assumptions in Theorem \ref{thm:nonlinear}, $H[f|\theta_M] < \infty$ and since $H[\theta(\tau)|\theta_M]$ takes values continuously in time, we may define $[\tau_0,\tau_\star]$ to be the largest connected time interval such that the following bootstrap hypothesis holds:
\begin{align}
\sup_{\tau \in (\tau_0,\tau_\star)}H[\theta(\tau)|\theta_M] & \leq 4H[f|\theta_M]. \label{ineq:bootstrapNL} 
\end{align}  
By propagation of regularity and continuity in time, $\tau_0 < \tau_\star < \tau_{\max}$ \cite{CalvezCarrillo06,BRB10}.
The essential component of the proof of Theorem \ref{thm:nonlinear} is to prove that $\tau_\star = \infty$. 
Ultimately, we will be able to choose $\tau_0$ large enough such that on $(\tau_0,\tau_\star)$, $H[\theta(\tau)|\theta_M] < 2H[f|\theta_M]$, and hence $\tau_\star = \infty$.   

The first step is to compute the time evolution of the relative entropy as for instance in 
\cite{CanizoCarrilloSchonbek,BedrossianIA10} (note that these computations can be justified on $[\tau_0,\tau_{\max})$ by propagation of regularity \cite{Blanchet09,BRB10}).
By Cauchy-Schwarz and the definition of the entropy production functional $I$ \eqref{def:I}, we have the following:
\begin{align}
\frac{d}{d\tau}H[\theta(\tau)|\theta_M] & = -I[\theta] + e^{(N-1)\tau} \int \grad \left(\frac{m\theta^{m-1}}{m-1} + \frac{1}{2}\abs{\xi}^2\right) \cdot \theta \grad \Kk(e^{\tau} \cdot) \ast \theta d\xi \nonumber \\ 
& \leq -I[\theta] + e^{(N-1)\tau}I[\theta]^{1/2}\sqrt{\int \theta \abs{\grad \Kk(e^\tau \cdot) \ast \theta}^2 d\xi}. \label{ineq:Hevo}
\end{align}  
The latter term is an error that we must control in order to propagate \eqref{ineq:bootstrapNL}. 
By H\"older's inequality and Young's inequality:   
\begin{align} 
\sqrt{\int \theta \abs{\grad \Kk(e^\tau \cdot) \ast \theta}^2 d\xi} \leq \norm{\theta}_m^{1/2}\norm{\grad \Kk(e^\tau \cdot) \ast \theta}_{\frac{2m}{m-1}} \lesssim e^{-\frac{d\tau}{q}}\norm{\grad \Kk}_q \norm{\theta}_m^{1/2} \norm{\theta}_{p}, \label{ineq:symbreakNL}
\end{align}
where here $p \in \left[ \frac{2md}{md + 2m - d}, \frac{2m}{m-1}\right)$ satisfies
\begin{align} 
\frac{1}{p} = 1 + \frac{m-1}{2m} - \frac{1}{q}. \label{def:p}
\end{align}
Note that if $q = 1$, then $p = \frac{2m}{m-1}$; also note that for no choice of $d \geq 3$ do we get $p \leq m$ (since $m = 2-2/d$).  
Applying \eqref{ineq:symbreakNL} to the evolution of the relative entropy \eqref{ineq:Hevo} implies that for some constant $C> 0$ depending on $\Kk$, 
\begin{align*} 
\frac{d}{d\tau}H[\theta(\tau)|\theta_M] & \leq -I[\theta] + C e^{\left(d-1 - \frac{d}{q}\right)\tau}I[\theta]^{1/2} \norm{\theta}_m^{1/2} \norm{\theta}_{p}. 
\end{align*} 
The exponent is negative due to the assumption that $q < \frac{d}{d-1}$ and this will provide
 the small parameter which we may use to close the bootstrap argument.  
For notational simplicity denote
\begin{align*} 
\epsilon = -\left(d-1 - \frac{d}{q}\right) > 0. 
\end{align*}
Since, 
\begin{align*} 
\frac{1}{m-1}\norm{\theta}_m^m \leq H[\theta|\theta_M] + H[\theta_M], 
\end{align*} 
we have (adjusting $C$ each line), 
\begin{align*} 
\frac{d}{d\tau}H[\theta(\tau)|\theta_M] & \leq -I[\theta] + Ce^{-\epsilon\tau}I[\theta]^{1/2} \left(H[\theta|\theta_M]^{\frac{1}{2m}} + H[\theta_M]^{\frac{1}{2m}} \right) \norm{\theta}_{p} \\ 
& \leq -\frac{1}{2}I[\theta]  + Ce^{-2\epsilon\tau}\left(H[\theta|\theta_M]^{\frac{1}{m}} + H[\theta_M]^{\frac{1}{m}} \right) \norm{\theta}^2_{p} \\ 
& \leq -\frac{1}{2}I[\theta] + \frac{1}{4}H[\theta|\theta_M]  + Ce^{-\frac{2m}{(m-1)}\epsilon\tau}\norm{\theta}^{\frac{2m}{m-1}}_{p} + C H[\theta_M]^{\frac{1}{m}}e^{-2\epsilon\tau} \norm{\theta}_{p}^2. 
\end{align*}
Applying the crucial \eqref{ineq:relative_entropy} then implies
\begin{align*} 
\frac{d}{d\tau}H[\theta(\tau)|\theta_M] &  \leq -\frac{3}{4}H[\theta|\theta_M] + C e^{-\frac{2m}{(m-1)}\epsilon\tau}\norm{\theta}^{\frac{2m}{m-1}}_{p} + C H[\theta_M]^{\frac{1}{m}}e^{-2\epsilon\tau} \norm{\theta}_{p}^2.
\end{align*}  
Integrating this over $(\tau_0,\tau_\star)$ gives (adjusting $C$ again)
\begin{align} 
\sup_{\tau \in (\tau_0,\tau_\star)}H[\theta(\tau)|\theta_M] & \leq H[f|\theta_M] + C e^{-\frac{2m}{(m-1)}\epsilon\tau_0}\left(\sup_{\tau \in (\tau_0,\tau_\star)}\norm{\theta(\tau)}^{\frac{2m}{m-1}}_{p}\right) \nonumber\\ & \quad 
+ C e^{-2\epsilon\tau_0} \left(\sup_{\tau \in (\tau_0,\tau_\star)} \norm{\theta(\tau)}_{p}^2 \right). \label{ineq:Hcontrol}
\end{align}

Since $p > m$, in order to control the RHS of \eqref{ineq:Hcontrol}, we need a second estimate on the high norm $L^p$.
This estimate will be obtained by truncated $L^p$ estimate methods; we will especially model the arguments after those found in \cite{Kowalczyk05,Blanchet09,CalvezCarrillo06,BRB10}.   
The necessary equi-integrability will come from \eqref{ineq:Hcontrol}, coupling the high and low norm estimates together.
Then $\tau_0$ will be chosen large in order to close the argument. 
 
Denote $\theta_k := (\theta - k)_+$ and recall that for all $1 \leq r < \infty$:  
\begin{equation}
\norm{\theta}_r^r \lesssim_r \norm{\theta_k}_r^r + k^{r-1}\norm{\theta}_1. \label{ineq:slicing}
\end{equation}
Compute the evolution of $\norm{\theta_k}_p^p$, using that $\theta^l_k \theta = \theta_k^{l+1} + k\theta_k^{l}$ and $\grad \theta^l = \grad \theta_k^l$ for all $l > 0$: 
\begin{align*} 
\frac{d}{d\tau}\norm{\theta_k(\tau)}_p^p & = -\frac{4mp(p-1)}{(p+m-1)^2}\int \abs{\nabla \theta_k^{\frac{p+m-1}{2}}}^2 d\xi - \int \left((p-1)\theta_k^{p} + kp\theta^{p-1}_k \right) \nabla \cdot \left(e^{(d-1)\tau}\nabla \K(e^{\tau}\cdot) \ast \theta \right)d\xi \\ 
& \quad + d(p-1)\norm{\theta_k}^{p+1}_{p+1} + dkp\norm{\theta_k}_{p}^{p}.
\end{align*} 
By H\"older's inequality, the Calderon-Zygmund inequality \cite{BigStein} (applied to the singular integral operator $e^{d\tau}\Delta\Kk(e^\tau \cdot)$ -- one can verify that the constants do not depend on $\tau$ \cite{BedrossianIA10}) and \eqref{ineq:slicing} (again adjusting $C$  every line): 
\begin{align*} 
\frac{d}{d\tau}\norm{\theta_k(\tau)}_p^p & \leq -\frac{4mp(p-1)}{(p+m-1)^2}\int \abs{\nabla \theta_k^{\frac{p+m-1}{2}}}^2 d\xi 
+ (p-1)\norm{\theta_k}_{p+1}^p\norm{e^{d\tau} \Delta\K(e^\tau \cdot) \ast \theta}_{p+1} \\ & \quad + kp\norm{\theta_k}_{p}^{p-1}\norm{e^{d\tau} \Delta \K(e^\tau \cdot) \ast \theta}_{p} + d(p-1)\norm{\theta_k}_{p+1}^{p+1} + dkp\norm{\theta_k}_{p}^{p} \\  
& \leq -\frac{4mp(p-1)}{(p+m-1)^2}\int \abs{\nabla \theta_k^{\frac{p+m-1}{2}}}^2 d\xi 
 + C(p,d,\K)\norm{\theta_k}_{p+1}^{p+1} 
 + C(p,d,k,\K)\norm{\theta_k}_{p}^p \\
& \leq -\frac{4mp(p-1)}{(p+m-1)^2}\int \abs{\nabla \theta_k^{\frac{p+m-1}{2}}}^2 d\xi + C_A\norm{\theta_k}_{p+1}^{p+1} + C_L, 
\end{align*}
where the last line followed by interpolation and we are defining the constants $C_A$ (which depends on $\Kk$, $d$ and $p$) 
and $C_L$ (which depends on $d,k,M,\K$ and $p$) for future convenience.  
By an appropriate Gagliardo-Nirenberg-Sobolev inequality, as in \cite{CalvezCarrillo06,Blanchet09,BRB10,BedrossianIA10}, 
we have for some constant $C_D$ (depending ultimately on $d$ and $p$),
\begin{align} 
\frac{d}{d\tau}\norm{\theta_k(\tau)}_p^p \leq  \left(-\frac{C_D}{\norm{\theta_k}_1^{2-m}} + C_A\right)\norm{\theta_k}_{p+1}^{p+1} + C_L.  \label{ineq:Lpctrl}
\end{align}  
The key point here is that control on $H[\theta|\theta_M]$ implies that $\norm{\theta_k}_1$ will decrease at a known rate with increasing $k$ (equivalent to equi-integrability) and hence used to make the first term a priori negative. Indeed,   
\begin{align} 
\norm{\theta_k}_1 \leq k^{1-m}\norm{\theta}_m^m \lesssim k^{1-m}\left(H[\theta(\tau)|\theta_M]) + H[\theta_M]\right). \label{ineq:equicontrol}
\end{align}
Therefore, by \eqref{ineq:bootstrapNL}, we can pick a $k=k_0(H[f|\theta_M],M)$ sufficiently large depending only on $d$,$H[f|\theta_M]$, $M$, $p$ and $\Kk$ (via $C_A$) such that on $(\tau_0,\tau_\star)$ we have 
\begin{align*} 
-\frac{C_D}{\norm{\theta_k}_1^{2-m}} + C_A < -1. 
\end{align*}
Hence by \eqref{ineq:Lpctrl} and the interpolation $\norm{\theta}_p^p \leq \norm{\theta}_{p+1}^{p+1} + M$ (note $C_L$ is now fixed large depending on $k_0$) 
\begin{align*} 
\frac{d}{d\tau}\norm{\theta_k(\tau)}_p^p & \leq -\norm{\theta_k}_{p+1}^{p+1} + C_L \\ 
&\leq -\norm{\theta_k}_{p}^{p} + M + C_L.
\end{align*}
Upon integration, this yields the following: 
\begin{align*} 
\sup_{\tau \in (\tau_0, \tau_\star)} \norm{\theta_k(\tau)}_p^p \leq \max\left( \norm{f_k}_p^p, M + C_L\right). 
\end{align*} 
By \eqref{ineq:slicing} it follows that 
\begin{align} 
\sup_{\tau \in (\tau_0, \tau_\star)} \norm{\theta(\tau)}_p^p \lesssim_p \max\left( \norm{f_k}_p^p, M + C_L\right) + k_0^{p-1}M. \label{ineq:highnrm}
\end{align} 
Note that the constants do not depend on $\tau_\star$.
Applying the control \eqref{ineq:highnrm} in \eqref{ineq:Hcontrol} implies that over the time interval $[\tau_0,\tau_\star)$, for some $C_F = C_F(\norm{f}_p,H[f|\theta_M],M,\K,d,p)$, we have
\begin{align*} 
\sup_{\tau \in (\tau_0,\tau_\star)}H[\theta(\tau)|\theta_M] \leq H[f|\theta_M] + C_F e^{-\frac{2m}{(m-1)}\epsilon\tau_0}. 
\end{align*}
It follows that if we choose $\tau_0$ depending only on $C_F$ and $H[f|\theta_M]$ then,  
\begin{align} 
\sup_{\tau \in (\tau_0,\tau_\star)}H[\theta(\tau)|\theta_M] \leq 2H[f|\theta_M]. \label{ineq:glbHbd}
\end{align}
Hence $\tau_\star = \tau_{\max}$, which implies also that \eqref{ineq:highnrm} holds until $\tau_{\max}$. 
By the regularity theory for \eqref{def:ADD} it follows that $\tau_{\max} = \infty$ (see e.g. \cite{CalvezCarrillo06,BRB10}) and therefore both \eqref{ineq:highnrm} and \eqref{ineq:glbHbd} hold globally in time. 

Since \eqref{ineq:glbHbd} controls a norm with regularity higher than $L^1$ in the similarity variables \eqref{def:NLsim}, Theorem 1(ii) of \cite{BedrossianIA10} implies the optimal $L^\infty$ decay estimate \eqref{ineq:LinftyNL}. 
Theorems 2 or 3 of \cite{BedrossianIA10} further imply as well the convergence to the Barenblatt solution at the specific rate depending on the decay of the interaction potential as stated in \eqref{ineq:BBlattNL}. 
\end{proof} 

\section*{Acknowledgments}
The author would like to thank Adrien Blanchet, Jose A. Carrillo and Marco Di Francesco for helpful discussions. Partially supported by NSF Postdoctoral Fellowship in Mathematical Sciences, DMS-1103765.

\bibliographystyle{plain}
\bibliography{nonlocal_eqns}

\end{document}